\documentclass{article}
\usepackage{amsfonts}
\usepackage{amsmath}

\newtheorem{theorem}{Theorem}

\newtheorem{corollary}[theorem]{Corollary}

\newtheorem{lemma}[theorem]{Lemma}

\begin{document}

\title{Construction of general forms of ordinary generating functions for more families of numbers and multiple variables polynomials}
\author{Yilmaz SIMSEK \\
Department of Mathematics, Faculty of Science \\
University of Akdeniz TR-07058 Antalya, Turkey\\
ysimsek@akdeniz.edu.tr}
\maketitle

\begin{abstract}
The aim of this paper is to construct general forms of ordinary generating
functions for special numbers and polynomials involving Fibonacci type
numbers and polynomials, Lucas numbers and polynomials, Chebyshev
polynomials, Sextet polynomials, Humbert-type numbers and polynomials, chain
and anti-chain polynomials, rank polynomials of the lattices, length of any
alphabet of words, partitions, and other graph polynomials. By applying the Euler
transform and the Lambert series to these generating functions, many new
identities and relations are derived. By using differential equations of
these generating functions, some new recurrence relations for these
polynomials are found. Moreover, general Binet's type formulas for these polynomials
are given. Finally, some new classes of polynomials and their corresponding
certain family of special numbers are investigated with the help of these
generating functions.


\end{abstract}

\section{Introduction}

Special polynomials and numbers with their generating functions play a vital
role not only in many branches of mathematics, but also in almost many
applied sciences, thanks to their important applications. Therefore, the development of these functions is always destined to be up to date. As a result, many researchers have been doing very deep research on these functions for more than a century in order to  solve many different problems involving  mathematical models. Modified generating functions for the Fibonacci type polynomials,  the Lucas type polynomials, and partitions have been given by many researchers for almost regularly every year. Most of these studies consists of modifying or unifying existing generating functions by adding either parameters or a few polynomials to the coefficients of existing generating functions \cite{alkanMTJp,Askey,Belbeshir,Bravo,AliB,Ceyvin,Charalambos,comtet,Deza,djordjevic,Dj-Mi2014,Forgacs0,Forgacs1,Forgacs2,gegenbauer,Goubi,GoubiMTJPAM,Gould,gould,Horadam,Hosonay,humbert,KahanBMS,KahanVit,KilarSimsekASCM,koshy1,KucukoAADm,KucukogluTJM,KucukogluAxio,Lia,lother,Munarini,Ortac,GulsahFilomat,GSG,Rabolovic,Pinch,Simsek,Simsekaadm,Simsek.mmas,SrivatavaChoi,srivastava,SrivatavaJNT,Suetin,Uygun}.

For these fundamental reasons, the main motivation of this paper is to give more general generating functions for sequences of polynomials and numbers. Thus, we define the following generating functions for sequences of new classes of multiple variables polynomials, which denoted respectively by  $\mathbb{Y}_{n}\left( P(%
\overrightarrow{X_{m}})\right) $ and  $\mathbb{S}_{n}\left( P(\overrightarrow{%
	X_{m}});Q(\overrightarrow{X_{k}})\right) $:
\begin{equation}
	F\left( w,P(\overrightarrow{X_{m}})\right) =\frac{1}{1+\sum%
		\limits_{j=1}^{m}P_{j}(x_{j})w^{j}}=\sum_{n=0}^{\infty }\mathbb{Y}_{n}\left(
	P(\overrightarrow{X_{m}})\right) w^{n},  \label{Ag}
\end{equation}%
and%
\begin{eqnarray}
	G\left( w,P(\overrightarrow{X_{m}});Q(\overrightarrow{X_{k}})\right) =
	\frac{\sum\limits_{j=0}^{k}Q_{j}(x_{j})w^{j}}{1+\sum%
		\limits_{j=1}^{m}P_{j}(x_{j})w^{j}}  
	=\sum_{n=0}^{\infty }\mathbb{S}_{n}\left( P(\overrightarrow{X_{m}});Q(%
	\overrightarrow{X_{k}})\right) w^{n},  \label{ag-1}
\end{eqnarray}%
where $P(\overrightarrow{X_{m}})=\left( P_{1}(x_{1}),P_{2}(x_{2}),\ldots
,P_{m}(x_{m})\right) $, \ $Q(\overrightarrow{X_{k}})=\left(
Q_{1}(x_{1}),Q_{2}(x_{2}),\ldots ,Q_{k}(x_{k})\right) $, 
\begin{equation*}
	P_{j}(x_{j})=\sum\limits_{v=0}^{d}a_{v}x_{j}^{v}, \qquad Q_{l}(x_{l})=\sum\limits_{v=0}^{c}b_{v}x_{l}^{v}
\end{equation*} are any polynomials in $x_{j}$, $x_{l}$
and $m\in \mathbb{N}$, $c,d,k\in \mathbb{N}_{0}$, $0\leq l\leq k$ and $0\leq j\leq m$.

By using (\ref{Ag}) and (\ref{ag-1}), we investigate many properties of the
polynomials $\mathbb{Y}_{n}\left( P(\overrightarrow{X_{m}})\right) $ and $%
\mathbb{S}_{n}\left( P(\overrightarrow{X_{m}});Q(\overrightarrow{X_{k}}%
)\right) $.

In the next sections, we investigate in detail how these two new classes of polynomials are reduced to some special numbers and polynomials , with some special cases.

\subsection{Preliminaries}

In order to present our results, we need to give some special classes of
polynomials and numbers with their generating functions.

The Fibonacci polynomials, which are a polynomial sequence, are defined by the
following ordinary generating function:%
\begin{equation*}
	\frac{t}{1-xt-t^{2}}=\sum_{n=0}^{\infty }F_{n}\left( x\right) t^{n},
\end{equation*}%
where $F_{0}\left( x\right) =0$, $F_{1}\left( x\right) =1$ and $F_{2}\left(
x\right) =x$. It is clear that the degree of the polynomials is $n-1$. The
polynomials $F_{n}\left( x\right) $\ can also be stated in terms of not
only well-known the Lucas polynomials, but also the Chebyshev polynomials of the
second kind $U_{n-1}\left( \frac{ix}{2}\right) =i^{n-1}F_{n}\left( x\right) 
$, $i$ is a imaginary unit $i^{2}=-1$. These polynomials can also be taken
into account as a generalization of the Fibonacci numbers $F_{n}$.
That is $F_{n}:=F_{n}\left( 1\right) $. The Pell numbers are obtained by
evaluating $F_{n}\left( 2\right) $ (\textit{cf}. \cite[p. 411]{koshy1}; see also \cite{KilarSimsekASCM}).

The other important class of the polynomials, which are generated in a
similar way from the Lucas numbers, are known as Lucas polynomials. These
polynomials are defined by the
following ordinary generating function:
\begin{equation*}
	\frac{2-xt}{1-xt-t^{2}}=\sum_{n=0}^{\infty }L_{n}\left( x\right) t^{n},
\end{equation*}%
where $L_{0}\left( x\right) =2$, $L_{1}\left( x\right) =x$ and $L_{2}\left(
x\right) =x^{2}+2$ (\textit{cf}. \cite[p. 26]{koshy1}). It is clear that
the degree of the polynomials is $n$. The polynomials $L_{n}\left( x\right) $
can also be stated in terms of the Chebyshev polynomials of the first kind $2T_{n}\left( -\frac{ix}{2}\right)=i^{n}L_{n}\left( x\right)
$. These polynomials can also be
taken into account as a generalization of the Lucas numbers $%
L_{n}$. That is $L_{n}:=L_{n}\left( 1\right) $ (\textit{cf}. \cite[p. 372]%
{koshy1}; see also \cite{KilarSimsekASCM}).

The Humbert polynomials $\left\{ \Pi _{n,m}^{\lambda }\right\}
_{n=0}^{\infty }$ were defined in $1921$ by Humbert \cite{humbert}. Their
generating function is given by
\begin{equation}
	\frac{1}{\left( 1-mxt+t^{m}\right) ^{\lambda }}=\sum_{n=0}^{\infty }\Pi
	_{n,m}^{\lambda }\left( x\right) t^{n}.  \label{HP}
\end{equation}

The Fibonacci type polynomials in two variables are defined by the
following ordinary generating function:
\begin{equation}
	H\left( t;x,y;k,m,n\right) =\sum_{j=0}^{\infty }\mathcal{G}_{j}\left(
	x,y;k,m,n\right) t^{j}=\frac{1}{1-x^{k}t-y^{m}t^{m+n}},  \label{GH}
\end{equation}%
where $k,m,n\in \mathbb{N}_{0}$ (\textit{cf}. \cite{GulsahFilomat}). An explicit formula for the polynomials $%
\mathcal{G}_{j}\left( x,y;k,m,n\right) $ is given by%
\begin{equation*}
	\mathcal{G}_{j}\left( x,y;k,m,n\right) =\sum_{c=0}^{\left[ \frac{j}{m+n}%
		\right] }\binom{j-c\left( m+n-1\right) }{c}y^{mc}x^{jk-mck-nck},
\end{equation*}%
where $\left[ a\right] $ is the largest integer $\leq a$ (\textit{cf}. \cite{GulsahFilomat}).

Two variable Fibonacci type polynomials of higher order are defined by the
following generating function 
\begin{equation}
	\sum_{j=0}^{\infty }\mathcal{G}_{j}^{\left( h\right) }\left(
	x,y;k,m,n\right) t^{j}=\frac{1}{\left( 1-x^{k}t-y^{m}t^{n+m}\right) ^{h}},
	\label{r2}
\end{equation}%
where $h$ is a positive integer (\textit{cf.} \cite{GSG}). Observe that%
\begin{eqnarray*}
	\mathcal{G}_{j}^{\left( 1\right) }\left( x,y;k,m,n\right) &=&\mathcal{G}%
	_{j}\left( x,y;k,m,n\right) \\
	\mathcal{G}_{j}^{\left( h\right) }\left( ax,-1;1,1,a-1\right) &=&\Pi
	_{j,a}^{h}\left( x\right)
\end{eqnarray*}%
(\textit{cf.} \cite{GSG}).

The results of this paper are briefly organized as follows. 

In Section \ref{section2}, by using
the Eqs. (\ref{Ag}) and (\ref{ag-1}), we give explicit formulas for the polynomials $\mathbb{Y}%
_{n}\left( P(\overrightarrow{X_{m}})\right) $ and $\mathbb{S}_{n}\left( P(%
\overrightarrow{X_{m}});Q(\overrightarrow{X_{k}})\right)$. We also give some relations and identities for these polynomials.
In Section \ref{section3}, we show that many certain families of special numbers and polynomials can be given in
terms of the polynomials $\mathbb{Y}_{n}\left( P(\overrightarrow{X_{m}}%
)\right) $ and $\mathbb{S}_{n}\left( P(\overrightarrow{X_{m}});Q(%
\overrightarrow{X_{k}})\right)$. We also give some special values of these polynomials. Moreover, we give Binet type formulas for the polynomials. 
In Sections\ref{section4}, we construct generating functions for higher order of the polynomials $\mathbb{Y}%
_{n}\left( P(\overrightarrow{X_{m}})\right) $ and $\mathbb{S}_{n}\left( P(%
\overrightarrow{X_{m}});Q(\overrightarrow{X_{k}})\right) $. We give some
properties of these polynomials with their generating functions and their special values.
Finally, Section \ref{section5}, we give some partial derivative equations of the generating functions for the polynomials  $\mathbb{Y}%
_{n}\left( P(\overrightarrow{X_{m}})\right) $ and $\mathbb{S}_{n}\left( P(%
\overrightarrow{X_{m}});Q(\overrightarrow{X_{k}})\right)$. By using these equation, we derive some recurrence relations for these polynomials.

\section{Explicit formulas for polynomials  $\mathbb{Y}_{n}\left( P(\protect%
	\overrightarrow{X_{m}})\right) $ and $\mathbb{S}_{n}\left( P(\protect%
	\overrightarrow{X_{m}});Q(\protect\overrightarrow{X_{k}})\right) $} \label{section2}

In this section, we give some explicit formulas for the polynomials $\mathbb{Y}%
_{n}\left( P(\overrightarrow{X_{m}})\right) $ and $\mathbb{S}_{n}\left( P(%
\overrightarrow{X_{m}});Q(\overrightarrow{X_{k}})\right) $ with the aid of
the Eqs. (\ref{Ag}) and (\ref{ag-1}).

Under the convergence conditions of the geometric series, we give explicit
formulas for polynomials $\mathbb{Y}_{n}\left( P(\overrightarrow{X_{m}}%
)\right) $ and $\mathbb{S}_{n}\left( P(\overrightarrow{X_{m}});Q(%
\overrightarrow{X_{k}})\right) $.

Putting $m=2$ in (\ref{Ag}), we obtain%
\begin{equation}
	\frac{1}{1+P_{1}(x_{1})w+P_{2}(x_{2})w^{2}}=\sum_{n_{1}=0}^{\infty }\mathbb{Y%
	}_{n_{1}}\left( P(\overrightarrow{X_{2}})\right) w^{n_{1}}.  \label{m2}
\end{equation}%
Using (\ref{m2}), we have%
\begin{eqnarray*}
	\sum_{n_{1}=0}^{\infty }\mathbb{Y}_{n_{1}}\left( P(\overrightarrow{X_{2}}%
	)\right) w^{n_{1}}&=&\sum_{n_{1}=0}^{\infty }\sum_{n_{2}=0}^{\infty
	}(-1)^{n_{1}+n_{2}}\binom{n_{1}+n_{2}}{n_{2}}\left( P_{1}(x_{1})\right)
	^{n_{1}}  \\  &&\times \left( P_{2}(x_{2})\right) ^{n_{2}}w^{n_{1}+2n_{2}}.
\end{eqnarray*}%
Replacing respectively $n_{1}$ and $n_{2}$ by $n_{1}-2n_{2}\ $\ and $n_{2}$
in the interior of the above sums, we get 
\begin{eqnarray*}
	\sum_{n_{1}=0}^{\infty }\mathbb{Y}_{n_{1}}\left( P(\overrightarrow{X_{2}}%
	)\right) w^{n_{1}} &=&\sum_{n_{1}=0}^{\infty }\sum_{n_{2}=0}^{\left[ \frac{%
			n_{1}}{2}\right] }(-1)^{n_{1}-n_{2}}\binom{n_{1}-n_{2}}{n_{2}}\left(
	P_{1}(x_{1})\right) ^{n_{1}-2n_{2}} \\
	&&\times \left( P_{2}(x_{2})\right) ^{n_{2}}w^{n_{1}}.
\end{eqnarray*}%
Comparing the coefficients of $w^{n_{1}}$ on both sides of the above
equation yields the following formula for the polynomials $\mathbb{Y}%
_{n_{1}}\left( P_{1}(x),P_{2}(x_{2})\right) $:
\begin{lemma}\label{lemma1}
	Let  $n_{1}\in \mathbb{N}_{0}$. Then we have%
	\begin{eqnarray}
		\mathbb{Y}_{n_{1}}\left( P_{1}(x),P_{2}(x_{2})\right) &=&\sum_{n_{2}=0}^{\left[
			\frac{n_{1}}{2}\right] }(-1)^{n_{1}-n_{2}}\binom{n_{1}-n_{2}}{n_{2}}\left(
		P_{1}(x_{1})\right) ^{n_{1}-2n_{2}}\label{m2a}\\ &&\times \left( P_{2}(x_{2})\right) ^{n_{2}}.\notag
	\end{eqnarray}
\end{lemma}
Substituting $m=3$ into (\ref{Ag}), we get%
\begin{equation*}
	\frac{1}{1+P_{1}(x_{1})w+P_{2}(x_{2})w^{2}+P_{3}(x_{3})w^{3}}%
	=\sum_{n_{1}=0}^{\infty }\mathbb{Y}_{n_{1}}\left( P(\overrightarrow{X_{3}}%
	)\right) w^{n_{1}}.
\end{equation*}%
After some calculations in the above generating function, we obtain%
\begin{eqnarray*}
	&&\sum_{n_{1}=0}^{\infty }\mathbb{Y}_{n_{1}}\left( P(\overrightarrow{X_{3}}%
	)\right) w^{n_{1}} \\
	&=&\sum_{n_{1}=0}^{\infty }\sum_{n_{2}=0}^{\infty }\sum_{n_{3}=0}^{\infty }(-1)^{n_{1}+n_{2}+n_{3}}%
	\binom{n_{1}+n_{2}}{n_{2}}\binom{n_{1}+n_{2}+n_{3}}{n_{3}}\\
	&&\times
	\prod\limits_{j=1}^{3}\left( P_{j}(x_{j})\right)
	^{n_{j}}w^{n_{1}+2n_{2}+3n_{3}}.
\end{eqnarray*}%
Replacing respectively $n_{1}$, $n_{2}$ and $n_{3}$ by $%
n_{1}-2(n_{2}-2n_{3})-3n_{3}$, $n_{2}-2n_{3}\ $\ and $n_{3}$ in the interior
of the above sums, we get 
\begin{eqnarray*}
	&&\sum_{n_{1}=0}^{\infty }\mathbb{Y}_{n_{1}}\left( P(\overrightarrow{X_{3}}%
	)\right) w^{n_{1}} \\
	&=&\sum_{n_{1}=0}^{\infty }\sum_{n_{2}=0}^{\left[ \frac{n_{1}}{2}\right]
	}\sum_{n_{3}=0}^{\left[ \frac{n_{2}}{3}\right] }(-1)^{n_{1}-n_{2}}\binom{%
		n_{1}-n_{2}-n_{3}}{n_{2}-2n_{3}}\binom{n_{1}-n_{2}}{n_{3}}\left(
	P_{1}(x_{1})\right) ^{n_{1}-2n_{2}+n_{3}} \\
	&&\times \left( P_{2}(x_{2})\right) ^{n_{2}-2n_{3}}\left(
	P_{3}(x_{3})\right) ^{n_{3}}w^{n_{1}}.
\end{eqnarray*}%
Comparing the coefficients of $w^{n_{1}}$ on both sides of the above
equation yields the following formula for the polynomials $\mathbb{Y}%
_{n_{1}}\left( P_{1}(x),P_{2}(x_{2}),P_{3}(x_{3})\right) $:
\begin{lemma} \label{lemma2}
	Let  $n_{1}\in \mathbb{N}_{0}$. Then we have%
	\begin{eqnarray}
		\mathbb{Y}_{n_{1}}\left( P(\overrightarrow{X_{3}})\right) &=&\sum_{n_{2}=0}^{%
			\left[ \frac{n_{1}}{2}\right] }\sum_{n_{3}=0}^{\left[ \frac{n_{2}}{3}\right]
		}(-1)^{n_{1}-n_{2}}\binom{n_{1}-n_{2}-n_{3}}{n_{2}-2n_{3}}\binom{n_{1}-n_{2}%
		}{n_{3}}  \label{m2ab} \\
		&&\times \left( P_{1}(x_{1})\right) ^{n_{1}-2n_{2}+n_{3}}\left( P_{2}(x_{2})\right)
		^{n_{2}-2n_{3}}\left( P_{3}(x_{3})\right) ^{n_{3}}.  \notag
	\end{eqnarray}
\end{lemma}
For $m\geq 4$, by using (\ref{Ag}), Lemma \ref{lemma1} and Lemma \ref{lemma2} with mathematical induction, we obtain%
\begin{equation*}
	\sum_{n_{1}=0}^{\infty }\frac{(-1)^{n_{1}}\left( P_{1}(x_{1})w\right)
		^{n_{1}}}{\left( 1+\sum\limits_{j=2}^{m}P_{j}(x_{j})w^{j}\right) ^{1+n_{1}}}%
	=\sum_{n_{1}=0}^{\infty }\mathbb{Y}_{n_{1}}\left( P(\overrightarrow{X_{m}}%
	)\right) w^{n_{1}}.
\end{equation*}%
From the above equation, we obtain%
\begin{eqnarray*}
	&&\sum_{n_{1}=0}^{\infty }\sum_{n_{2}=0}^{\infty }(-1)^{n_{1}+n_{2}}\binom{%
		n_{1}+n_{2}}{n_{2}}\frac{\left( P_{1}(x_{1})\right) ^{n_{1}}\left(
		P_{2}(x_{2})\right) ^{n_{2}}w^{n_{1}+2n_{2}}}{\left(
		1+\sum\limits_{j=3}^{m}P_{j}(x_{j})w^{j}\right) ^{1+n_{1}+n_{2}}} \\
	&=&\sum_{n_{1}=0}^{\infty }\mathbb{Y}_{n_{1}}\left( \underset{m\text{ times}}%
	{\underbrace{P_{1}(x_{1}),P_{2}(x_{2}),\ldots ,P_{m}(x_{m})}}\right)
	w^{n_{1}}.
\end{eqnarray*}%
Continuing these same processes sequentially, we arrive at the following
result:%
\begin{eqnarray*}
	&&\sum_{n=0}^{\infty }\mathbb{Y}_{n_{1}}\left( P(\overrightarrow{X_{m}}%
	)\right) w^{n_{1}} \\
	&&=\underset{m\text{ times}}{\underbrace{\sum_{n_{1}=0}^{\infty
			}\sum_{n_{2}=0}^{\infty }\cdots \sum_{n_{m}=0}^{\infty }}}\left(
	(-1)^{\sum_{j=1}^{m}n_{j}}\right) \binom{n_{1}+n_{2}}{n_{2}}\binom{%
		n_{1}+n_{2}+n_{3}}{n_{3}} \\
	&&\times \cdots \binom{n_{1}+n_{2}+\cdots +n_{m}}{n_{m}} \prod\limits_{v=1}^{m}\left( P_{v}(x_{v})\right)
	^{n_{v}}w^{\sum_{j=1}^{m}jn_{j}}.
\end{eqnarray*}

Replacing respectively $n_{1}$,$n_{2}$,\ldots ,$n_{m-1}$, and $n_{m}$ by $%
n_{1}-2(n_{2}-2n_{3})-3(n_{3}-3n_{4})-\cdots -\left( m-1\right)
(n_{m-1}-\left( m-1\right) n_{m})-mn_{m}$, $n_{2}-2n_{3}$, \ldots , $%
n_{m-1}-mn_{m}$, and $n_{m}$ in the interior of the above sums, and using the following known formula
\begin{equation*}
	\sum_{n_{1}=0}^{\infty }\sum_{n_{2}=0}^{\infty
	}A(n_{1},n_{2})=\sum_{n_{1}=0}^{\infty }\sum_{n_{2}=0}^{\left[ \frac{n_{1}}{a%
		}\right] }A(n_{1},n_{2}-an_{1}),
\end{equation*}%
where $a\in \mathbb{N}$, we obtain%
\begin{eqnarray*}
	&&\sum_{n_{1}=0}^{\infty }\mathbb{Y}_{n_{1}}\left( P(\overrightarrow{X_{m}}%
	)\right) w^{n_{1}}= \\
	&&\underset{m\text{ times}}{\underbrace{\sum_{n_{1}=0}^{\infty
			}\sum_{n_{2}=0}^{\left[ \frac{n_{1}}{2}\right] }\sum_{n_{3}=0}^{\left[ \frac{%
					n_{2}}{3}\right] }\cdots \sum_{n_{m}=0}^{\left[ \frac{n_{m-1}}{m}\right] }}}%
	\left( (-1)^{\sum_{j=1}^{m}n_{j}}\right) \binom{n_{1}+n_{2}}{n_{2}}\binom{%
		n_{1}+n_{2}+n_{3}}{n_{3}} \\
	&&\times\cdots \binom{n_{1}+n_{2}+\cdots +n_{m}}{n_{m}} \prod\limits_{v=1}^{m}\left( P_{v}(x_{v})\right) ^{n_{v}}w^{n_{1}}.
\end{eqnarray*}%
Comparing the coefficients of $w^{n_{1}}$ on both sides of the above
equation, we arrive at the following theorem:

\begin{theorem}
	Let $m\in \mathbb{N}$ with $m>1$ and $n_{1}\in \mathbb{N}_{0}$. Then we have%
	\begin{eqnarray}
		\mathbb{Y}_{n_{1}}\left( P(\overrightarrow{X_{m}})\right) =\left(
		\prod\limits_{j=2}^{m}\sum_{n_{j}=0}^{\left[ \frac{n_{j-1}}{j}\right]
		}\right) (-1)^{\sum_{j=1}^{m}n_{j}}\prod\limits_{d=2}^{m}\binom{\sum\limits_{k=1}^{d}n_{k}}{n_{d}}%
		\prod\limits_{v=1}^{m}\left( P_{v}(x_{v})\right) ^{n_{v}},\label{aT1}
	\end{eqnarray}%
	where%
	\begin{equation*}
		\left( \prod\limits_{j=2}^{m}\sum_{n_{j}=0}^{\left[ \frac{n_{j-1}}{j}\right]
		}\right) =\sum_{n_{2}=0}^{\left[ \frac{n_{1}}{2}\right] }\sum_{n_{3}=0}^{%
			\left[ \frac{n_{2}}{2}\right] }\cdots \sum_{n_{m}=0}^{\left[ \frac{n_{m-1}}{m%
			}\right] },
	\end{equation*}%
	
	\begin{equation*}
		\prod\limits_{d=2}^{m}\binom{\sum\limits_{k=1}^{d}n_{k}}{n_{d}}%
		=\binom{n_{1}+n_{2}}{n_{2}}\binom{%
			n_{1}+n_{2}+n_{3}}{n_{3}} \cdots \binom{n_{1}+n_{2}+\cdots +n_{m}}{n_{m}},
	\end{equation*}%
	
	\begin{equation*}
		\prod\limits_{v=1}^{m}\left( P_{v}(x_{v})\right) ^{n_{v}}=\left(
		P_{1}(x_{1})\right) ^{n_{1}}\left( P_{2}(x_{2})\right) ^{n_{2}}\cdots \left(
		P_{m}(x_{m})\right) ^{n_{m}},
	\end{equation*}
	\begin{equation*}
		(-1)^{\sum_{j=1}^{m}n_{j}}=(-1)^{n_{1}+n_{2}\cdots+n_{m}}.
	\end{equation*}%
	
\end{theorem}

Some special values of the Eqs. (\ref{Ag}) and (\ref{aT1}) are given as
follows:

Substituting $m=1$ into Eq. (\ref{Ag}), we get%
\begin{equation*}
	\frac{1}{1+P_{1}(x_{1})w}=\sum_{n=0}^{\infty }\mathbb{Y}_{n}\left(
	P_{1}(x_{1})\right) w^{n}.
\end{equation*}%
Therefore%
\begin{equation*}
	\mathbb{Y}_{n}\left( P_{1}(x_{1})\right) =\left( -P_{1}(x_{1})\right) ^{n}.
\end{equation*}

Substituting $P_{1}(x)=x_{1}$, $P_{2}(x)=x_{2}$, $P_{3}(x)=x_{3}$,\ldots ,$%
P_{m}(x)=x_{m}$ into Eq. (\ref{Ag}), we obtain%
\begin{equation*}
	\frac{1}{1+x_{1}w+x_{2}w^{2}+x_{3}w^{3}+\cdots +x_{m}w^{m}}%
	=\sum_{n_{1}=0}^{\infty }\mathbb{Y}_{n_{1}}\left( x_{1},x_{2},x_{3},\ldots
	,x_{m}\right) w^{n_{1}}.
\end{equation*}%
Combining the above equation with (\ref{aT1}), we obtain the following
result:

\begin{corollary}
	Let $m\in \mathbb{N}$ with $m>1$  and $n_{1}\in \mathbb{N}_{0}$. Then we have%
	\begin{equation}
		\mathbb{Y}_{n_{1}}\left( x_{1},x_{2},x_{3},\ldots ,x_{m}\right) =\left(
		\prod\limits_{j=2}^{m}\sum_{n_{j}=0}^{\left[ \frac{n_{j-1}}{j}\right]
		}\right) (-1)^{\sum_{j=1}^{m}n_{j}}\prod\limits_{d=2}^{m}\binom{%
			\sum\limits_{k=1}^{d}n_{k}}{n_{d}}\prod\limits_{j=1}^{m}x_{j}^{n_{j}}.
		\label{aT2}
	\end{equation}
\end{corollary}

Substituting $x_{1}=x$, $x_{2}=x^{2}$, $x_{3}=x^{3}$,\ldots , $x_{m}=x^{m}$
into (\ref{aT2}), we get the following result:

\begin{corollary}
	Let $m\in \mathbb{N}$ with $m>1$ and $n_{1}\in \mathbb{N}_{0}$. Then we have%
	\begin{eqnarray}
		\mathbb{Y}_{n_{1}}\left( x,x^{2},x^{3},\ldots ,x^{m}\right) &=&\left(
		\prod\limits_{j=2}^{m}\sum_{n_{j}=0}^{\left[ \frac{n_{j-1}}{j}\right]
		}\right) (-1)^{\sum_{j=1}^{m}n_{j}}	\label{ac1} \\ &&\times \prod\limits_{d=2}^{m}\binom{%
			\sum\limits_{k=1}^{d}n_{k}}{n_{d}} x^{n_{1}+2n_{2}+3n_{3}+\cdots +mn_{m}}. \notag
	\end{eqnarray}
\end{corollary}

Substituting $x_{1}=x$, $x_{2}=x$, $x_{3}=x$,\ldots , $x_{m}=x$ into (\ref%
{aT2}), we get the following result:

\begin{corollary}
	Let $m\in \mathbb{N}$ with $m>1$ and $n_{1}\in \mathbb{N}_{0}$. Then we have%
	\begin{eqnarray}
		\mathbb{Y}_{n_{1}}\left( \underset{m\textup{ times}}{\underbrace{x,x,\ldots ,x}%
		}\right) &=&\left( \prod\limits_{j=2}^{m}\sum_{n_{j}=0}^{\left[ \frac{n_{j-1}}{%
				j}\right] }\right) (-1)^{\sum_{j=1}^{m}n_{j}}	\label{acc11} \\ &&\times \prod\limits_{d=2}^{m}\binom{%
			\sum\limits_{k=1}^{d}n_{k}}{n_{d}}x^{n_{1}+n_{2}+n_{3}+\cdots +n_{m}}. \notag
	\end{eqnarray}
\end{corollary}

Substituting $x=1$ into (\ref{ac1}), we arrive at the following new general
family of special numbers or so-called certain classes of finite sum:

\begin{corollary}
	Let $n_{1}\in \mathbb{N}_{0}$. Then we have%
	\begin{equation}
		\mathbb{Y}_{n_{1}}\left( \underset{m\textup{ times}}{\underbrace{1,1,\ldots ,1}%
		}\right) =\left( \prod\limits_{j=2}^{m}\sum_{n_{j}=0}^{\left[ \frac{n_{j-1}}{%
				j}\right] }\right) (-1)^{\sum_{j=1}^{m}n_{j}}\prod\limits_{d=2}^{m}\binom{%
			\sum\limits_{k=1}^{d}n_{k}}{n_{d}}.  \label{ac2}
	\end{equation}
\end{corollary}

Substituting $x=-1$ into (\ref{ac1}), we also arrive at the following new
general family of special numbers or so-called certain classes of finite sum:

\begin{corollary}
	Let $n_{1}\in \mathbb{N}_{0}$. Then we have%
	\begin{eqnarray}
		\mathbb{Y}_{n_{1}}\left( \underset{m\textup{ times}}{\underbrace{-1,-1,\ldots
				,-1}}\right) &=&\left( \prod\limits_{j=2}^{m}\sum_{n_{j}=0}^{\left[ \frac{%
				n_{j-1}}{j}\right] }\right) (-1)^{\left(
			\sum_{j=1}^{m}n_{j}+n_{1}+2n_{2}+3n_{3}+\cdots +mn_{m}\right)
		} \notag	\\  &&\times \prod\limits_{d=2}^{m}\binom{\sum\limits_{k=1}^{d}n_{k}}{n_{d}}.
		\label{ac3}	\end{eqnarray}
\end{corollary}

Observe that the formula given in equation (\ref{ac3}) gives us the solution
of Exercise 21 in Charalambides's book \cite[p. 269]{Charalambos}, given
below:

Substituting $P_{j}(x_{j})=-1$, $j\in \left\{ 1,2,\ldots ,m\right\} $ into (%
\ref{Ag}), we get generating function for the Fibonacci numbers of order $m$:%
\begin{equation*}
	\frac{1}{1-w-w^{2}-\cdots -w^{m}}=\sum_{n=0}^{\infty }\mathbb{Y}_{n}\left( 
	\underset{m\text{ times}}{\underbrace{-1,-1,\ldots ,-1}}\right) w^{n},
\end{equation*}%
where%
\begin{equation*}
	F_{n,m}:=\mathbb{Y}_{n}\left( \underset{m\text{ times}}{\underbrace{%
			-1,-1,\ldots ,-1}}\right),
\end{equation*}%
where the number $F_{n+1,m}$ is the number of $n$-permutations of the set $\left\{ 0,1\right\} $ with
repetition and the restriction that no $m$ zeros are consecutive. By using (%
\ref{ac3}), we have  $F_{0,m}=F_{1,m}=1$. Hence, it is easy to show that%
\begin{equation*}
	F_{n,m}=2F_{n-1,m}-F_{n-m-1,m},
\end{equation*}%
where $n=m+1$, $m+2$, $m+3$, $\ldots $,%
\begin{equation*}
	F_{n,m}=\sum\limits_{v=0}^{\min \left\{ n,m\right\} }F_{n-v,m},
\end{equation*}%
(\textit{cf}. \cite[p. 269]{Charalambos}).

Recently, Tran \cite{Forgacs0} gave the following generating function (rational generating function) for sequence of polynomials $H_{v}(z)$:
\begin{eqnarray*}
	\sum_{v=0}^{\infty }H_{v}(z) w^{v}=\frac{1}{ 1+P_{1}(z)w+P_{m}(z)w^{m}},
\end{eqnarray*}%
where $P_{1}(z)$ and $P_{m}(z)$ are any polynomials in $z$ with complex coefficients. Combining (\ref{Ag}) with the above rational generating function, we have  $$H_{v}(z)=\mathbb{Y}_{v}\left( \underset{m-2\text{ times}}{P_{1}(z),\underbrace{0,0,\ldots ,0},P_{m}(z)}\right).$$ In \cite{Forgacs1} and \cite{Forgacs2}, Forgacs and Tran gave also many properties and applications of the rational generating functions.

By combining (\ref{ag-1}) with (\ref{Ag}), we get the following functional
equation:

\begin{equation*}
	G\left( w,P(\overrightarrow{X_{m}});Q(\overrightarrow{X_{k}})\right)
	=F\left( w,P(\overrightarrow{X_{m}})\right)
	\sum\limits_{j=0}^{k}Q_{j}(x_{j})w^{j}.
\end{equation*}%
Using the above equation, we get%
\begin{equation*}
	\sum_{n=0}^{\infty }\mathbb{S}_{n}\left( P(\overrightarrow{X_{m}});Q(%
	\overrightarrow{X_{k}})\right)
	w^{n}=\sum\limits_{j=0}^{k}Q_{j}(x_{j})\sum_{n=0}^{\infty }\mathbb{Y}%
	_{n}\left( P(\overrightarrow{X_{m}})\right) w^{n+j}.
\end{equation*}%
Therefore%
\begin{equation*}
	\sum_{n=0}^{\infty }\mathbb{S}_{n}\left( P(\overrightarrow{X_{m}});Q(%
	\overrightarrow{X_{k}})\right)
	w^{n}=\sum\limits_{j=0}^{k}Q_{j}(x_{j})\sum_{n=j}^{\infty }\mathbb{Y}%
	_{n-j}\left( P(\overrightarrow{X_{m}})\right) w^{n}.
\end{equation*}%
Comparing the coefficients of $w^{n}$ on both sides of the above
equation, we arrive at the following theorem:

\begin{theorem}
	Let $n,k\in \mathbb{N}_{0}$ with $n\geq k$ and $m\in \mathbb{N}$. Then we have%
	\begin{equation*}
		\mathbb{S}_{n}\left( P(\overrightarrow{X_{m}});Q(\overrightarrow{X_{k}}%
		)\right) =\sum\limits_{j=0}^{k}Q_{j}(x_{j})\mathbb{Y}_{n-j}\left( P(%
		\overrightarrow{X_{m}})\right) ,
	\end{equation*}%
	where $m$ and $k$ any nonnegative integers.
\end{theorem}

\section{Some special values of the polynomials $\mathbb{Y}_{n}\left( P(%
	\protect\overrightarrow{X_{m}})\right) $ and $\mathbb{S}_{n}\left( P(\protect%
	\overrightarrow{X_{m}});Q(\protect\overrightarrow{X_{k}})\right) $} \label{section3}

In this section, using some special values of polynomials $P(%
\overrightarrow{X_{m}})$ and $Q(\overrightarrow{X_{k}})$, we give special
values of the polynomials $\mathbb{Y}_{n}\left( P(\overrightarrow{X_{m}}%
)\right) $ and $\mathbb{S}_{n}\left( P(\overrightarrow{X_{m}});Q(%
\overrightarrow{X_{k}})\right) $ involving many certain families of the
special numbers and polynomials.

Some special certain families of numbers and polynomials can be expressed in
terms of the polynomials $\mathbb{Y}_{n}\left( P(\overrightarrow{X_{m}}%
)\right) $ and $\mathbb{S}_{n}\left( P(\overrightarrow{X_{m}});Q(%
\overrightarrow{X_{k}})\right)$. These numbers and polynomials are given as follows:

\textbf{Binomial coefficients}: Substituting $m=1$ and $%
P_{1}(x_{1})=-1-x_{1} $ into (\ref{Ag}), we have the following well-known
bivariate generating function for the binomial coefficients:%
\begin{eqnarray*}
	\frac{1}{1+\left( -1-x_{1}\right) w} =\sum_{n=0}^{\infty }\mathbb{Y}%
	_{n}\left( -1-x_{1}\right) w^{n} 
	=\sum_{n=0}^{\infty }\sum_{m=0}^{n }\binom{n}{m}x_{1}^{m}w^{n}.
\end{eqnarray*}
Therefore
\begin{eqnarray}
	\mathbb{Y}%
	_{n}\left( -1-x_{1}\right) = \sum_{m=0}^{n }\binom{n}{m}x_{1}^{m}. \label{Y22}
\end{eqnarray}

Substituting $x_{1}=1$ and $x_{1}=-1$ into (\ref{Y22}), respectively, we easily have $\mathbb{Y}%
_{n}\left( -2\right) = 2^{n}$ and $\mathbb{Y}%
_{n}\left( 0\right) = 0$.

\textbf{Sextet polynomials of hexagonal systems: }Substituting $%
P_{1}(x)=-x^{2}-4x-1,P_{2}(x_{2})=x^{2}$ into (\ref{m2}), we get the
following generating function for the sextet polynomials:%
\begin{equation}
	\frac{1}{1+(-x^{2}-4x-1)w+x^{2}w^{2}}=\sum_{n_{1}=0}^{\infty }\mathbb{Y}%
	_{n_{1}}\left( -x^{2}-4x-1,x^{2}\right) w^{n_{1}}.  \label{m2b}
\end{equation}%
Thus, combining the above equation with (\ref{m2a}), we obtain the following
explicit formula for the sextet polynomials:%
\begin{eqnarray}
	\mathbb{Y}_{n_{1}}\left( -x^{2}-4x-1,x^{2}\right) &=&\sum_{n_{2}=0}^{\left[ 
		\frac{n_{1}}{2}\right] }(-1)^{2n_{1}-3n_{2}}\binom{n_{1}-n_{2}}{n_{2}} \label{m2c}\\ &&\times 	\left(
	x^{2}+4x+1\right) ^{n_{1}-2n_{2}}x^{2n_{2}}.  \notag
\end{eqnarray}
Li et al. \cite{Lia} also gave generating function for the sextet
polynomials. They also gave some properties of these polynomials. The sextet
polynomials of hexagonal systems are symmetric, unimodal, log-concave, and
asymptotically normal.

The sextet polynomials are related to the Clar covering polynomial of
hexagonal systems and chromatic polynomials, which have many applications in
graphs theory and and in the theory of hexagonal systems. The sextet
polynomials are also very important to analyze hexagonal system or benzenoid
system with the aid of finite connected plane graph without cut vertices
with every interior face is bounded by a regular hexagon of side length $1$.
Moreover, topological properties of hexagonal systems have many important
applications in various quantum mechanical models of the electronic
structure of benzenoid hydrocarbons and also in resonance theory, Huckel
molecular orbital theory, Clar's aromatic sextet theory and the theory of
conjugated circuits (\textit{cf}. \cite{Ceyvin,Hosonay,Lia}).

\textbf{Rank polynomials: }Substituting $Q_{0}(x)=1$, $Q_{1}(x)=-x^{2}$ and $%
P_{1}(x)=-1-x-x^{2}$, $P_{2}(x)=x^{2}$, $P_{3}(x_{3})=x^{3}$ into (\ref{ag-1}%
), we get the following rank polynomials of the lattices $\mathcal{J}(%
\mathcal{G}_{m})$, which is the set of all ideals of the partially ordered
set $\mathcal{G}_{m}$ which is known as garland of order $m$:%
\begin{equation*}
	g_{m}(x)=\mathbb{S}_{m}\left( -1-x-x^{2},x^{2},x^{3};1,-x^{2}\right),
\end{equation*}%
where $	g_{m}(x)=\sum\limits_{k=0}^{m}g_{mk}x^{k}$. The rank polynomials of the lattices have various interesting applications
in the theory of Ockham algebras, and also in combinatorics involving fences
(or zigzag posets), crowns, garlands (or double fences) and many of their
generalizations of posets. The garland of order $m$ is the partially ordered set $\mathcal{G}_{m}$
briefly defined as follows: $\mathcal{G}_{0}$ is the empty poset, $\mathcal{G%
}_{1}$ is the chain of length $1$ and, for any other $m\geq 2$, $\mathcal{G}%
_{m}$ is the poset on $2m$ elements $x_{1}$, $x_{2}$,\ldots ,$x_{m}$ and $%
y_{1}$, $y_{2}$,\ldots ,$y_{m}$ with cover relations $x_{1}<y_{1}$, $%
x_{2}<y_{2}$,\ldots ,$x_{j}<y_{j-1}$, $x_{j}<y_{j}$ and $x_{j}<y_{j+1}$, for 
$j=2$, $3$,\ldots ,$m-1$ and $x_{m}<y_{m-1}$ and $x_{m}<y_{m}$. The ideals
and the anti-chains of a garland can be represented as words of a regular
language. Here the maximal elements of an ideal of a poset is an anti-chain,
and this establishes a bijective correspondence between all ideals of poset (%
\textit{cf}. \cite{Munarini}). An explicit formula for the number $g_{m}$ of
all anti-chains or ideals of $\mathcal{G}_{m}$ is given as follows:%
\begin{equation*}
	g_{m}=\mathbb{S}_{m}\left( -2,-1;1,1\right) =\frac{\left( 1+\sqrt{2}\right)
		^{m+1}+\left( 1-\sqrt{2}\right) ^{m+1}}{2}
\end{equation*}%
(\textit{cf}. \cite{Munarini}). 
Note that the formula of the number $g_{m}$ of
all anti-chains is also related to Exercise 11 in Charalambides's book \cite[p. 267]%
{Charalambos}. Due to the Exercise 11,  we also observe that
$g_{m}$ is the number of $m$-permutations of the set $\left\{
0,1,2\right\} $ with repetition and the restriction now two zeros and no two
ones are consecutive.

\textbf{Figurate numbers: }We give relations between special values of the
polynomials $\mathbb{S}_{n}\left( P(\overrightarrow{X_{m}});Q(%
\overrightarrow{X_{k}})\right) $ and the figurate numbers.

The $j$-gonal numbers, which are members of the space figurate numbers:%
\begin{equation*}
	\frac{x+(j-3)x^{2}}{1-3x+3x^{2}-x^{3}}=\sum\limits_{n=0}^{\infty }\mathbb{S}%
	_{n}\left( -3,3,-1;0,1,j-3\right) w^{n}
\end{equation*}%
is the generating functions for the sequence of hexagonal prism numbers.
That is,%
\begin{equation*}
	S_{j}(n):=\mathbb{S}_{n}\left( -3,3,-1;0,1,j-3\right) 
\end{equation*}%
(\textit{cf}. \cite[p. 145]{Deza}). For $j=3$, $S_{3}(n)$ reduces to the triangular numbers, which was studied
by Theon of Smyrna in the II-th century AC.

Hexagonal prism numbers, which are members of the space figurate numbers:%
\begin{equation*}
	\frac{x+10x^{2}+7x^{3}}{1-4x+6x^{2}-4x^{3}+x^{4}}=\sum\limits_{n=0}^{\infty }%
	\mathbb{S}_{n}\left( -4,6,-4,1;0,1,10,7\right) w^{n}
\end{equation*}%
is the generating function for the sequence of hexagonal prism numbers (%
\textit{cf}. \cite[p. 145]{Deza}). That is,%
\begin{equation*}
	PCS_{6}(n):=\mathbb{S}_{n}\left( -4,6,-4,1;0,1,10,7\right).
\end{equation*}%
The centered $j$-pyramidal numbers%
\begin{equation*}
	CS_{j}^{3}(n):=\mathbb{S}_{n}\left( -4,6,-4,1;0,1,j-2,1\right).
\end{equation*}%
The centered dodecahedron numbers%
\begin{equation*}
	\overline{D}(n):=\mathbb{S}_{n}\left( -4,6,-4,1;0,1,17,17,1\right).
\end{equation*}%
The centered icosahedron number%
\begin{equation*}
	\overline{I}(n):=\mathbb{S}_{n}\left( -4,6,-4,1;0,1,9,9,1\right).
\end{equation*}%
(Sloane's A005902).

The centered octahedron number%
\begin{equation*}
	\overline{O}(n):=\mathbb{S}_{n}\left( -4,6,-4,1;0,1,3,3,1\right),
\end{equation*}%
for detail, see \cite[p. 145]{Deza}.

\textbf{Anti-chain polynomials: }The anti-chain polynomials, which are also
related to not only the Chebyshev polynomials and the Sturm sequence but
also the Euler transform of a formal series, are given by%
\begin{equation*}
	G\left( w,-1-x,-x;1,x\right) =\sum\limits_{n=0}^{\infty }\mathbb{S}%
	_{n}\left( -1-x,-x;1,x\right) w^{n}.
\end{equation*}%
From the above equation, we get%
\begin{equation*}
	a_{n}(x)=\mathbb{S}_{n}\left( -1-x,-x;1,x\right) ,
\end{equation*}%
where%
\begin{equation*}
	a_{n}(x)=\sum\limits_{k=0}^{n}a_{n,k}x^{k}.
\end{equation*}%
These polynomials are log-concave and unimodal, for detail see \cite%
{Munarini}.

$\sum\limits_{n=0}^{\infty }a_{2n,n}x^{n}$ is the generating function,
associated with diagonal of the double series, for the following matrix,
which rows are given by the sequence A035607,%
\begin{equation*}
	A=\left[ a_{n,m}\right] =\left[ 
	\begin{array}{ccccccc}
		1 &  &  &  &  &  &  \\ 
		1 & 4 & 2 &  &  &  &  \\ 
		1 & 6 & 8 & 2 &  &  &  \\ 
		1 & 8 & 18 & 12 & 2 &  &  \\ 
		1 & 10 & 32 & 38 & 16 & 2 &  \\ 
		1 & 12 & 50 & 88 & 66 & 20 & 2 \\ 
		\ldots &  &  &  &  &  & 
	\end{array}%
	\right]
\end{equation*}%
is given by%
\begin{eqnarray*}
	G\left( w,-1-4x-x^{2},x^{2};1,x^{2}\right) &=&\sum\limits_{n=0}^{\infty }%
	\mathbb{S}_{n}\left( -1-4x-x^{2},x^{2};1,x^{2}\right) w^{n} \\
	&=&\frac{1+x^{2}w}{1+\left( -1-4x-x^{2}\right) w+x^{2}w^{2}},
\end{eqnarray*}%
for detail see \cite{Munarini}.

\textbf{Euler transform}: The Euler transform of a formal series $%
F(w)=\sum\limits_{j\geq 0}u_{j}w^{j}$ is given by%
\begin{eqnarray*}
	T^{\theta }\left( F(w)\right)=\frac{1}{1-\theta w}F\left( \frac{w}{%
		1-\theta w}\right) 
	=\sum\limits_{j=0}^{\infty }\sum\limits_{v=0}^{j}\binom{j}{v}\theta
	^{j-v}u_{v}w^{j}.
\end{eqnarray*}%
Thus by applying the Euler transform $T^{\theta }$ to the formal series $%
G\left( w,-1-x,-x;1,x\right) $ with respect to $w$, we obtain%
\begin{eqnarray}
	T^{x}\left\{ G\left( w,-1-x,-x;1,x\right) \right\} &=&F\left(
	w,-1-3x,2x^{2}\right)  \label{ac4} \\
	&=&F_{U}\left( \sqrt{2}xw,\frac{1+3x}{2\sqrt{2}x}\right) ,  \notag
\end{eqnarray}%
where $F_{U}\left( w,x\right) $ denotes the following generating function
for the Chebyshev polynomials of the second kind%
\begin{eqnarray*}
	F_{U}\left( w,x\right) =G\left( w,-2x,1;1\right) 
	=F\left( w,-2x,1\right)
	=\sum\limits_{n=0}^{\infty }U_{n}(x)w^{n},
\end{eqnarray*}%
where%
\begin{equation*}
	U_{n}(x):=\mathbb{S}_{n}\left( -2x,1;1\right) =\mathbb{Y}_{n}\left(
	-2x,1\right) .
\end{equation*}

Since the Euler transform is an invertible operator. By applying the Euler
invertible operator $T^{-\theta }=\left( T^{\theta }\right) ^{-1}$ to
Eq. (\ref{ac4}), we obtain%
\begin{eqnarray*}
	\mathbb{S}_{n}\left( -1-x,-x;1,x\right) =T^{-x}\left\{ F\left(
	w,-1-3x,2x^{2}\right) \right\} 
	=T^{-x}\left\{ F_{U}\left( \sqrt{2}xw,\frac{1+3x}{2\sqrt{2}x}\right)
	\right\} .
\end{eqnarray*}

For detail applications the Euler transform of a formal series to the
antichain polynomials and the Chebyshev polynomials, see \cite{Munarini}.

\textbf{Convolved Fermat quotients}: Goubi \cite{GoubiMTJPAM} gave
combinatorial formulation of the convolved Fermat quotients, which are
related to the polynomials $\mathbb{S}_{n}\left( P(\overrightarrow{X_{2m}}%
);Q(\overrightarrow{X_{k}})\right) $.

\textbf{Fibonacci polynomials}: The Fibonacci polynomials can also be
expressed in terms of the polynomials $\mathbb{Y}_{n}\left( P(%
\overrightarrow{X_{m}})\right) $ and $\mathbb{S}_{n}\left( P(\overrightarrow{%
	X_{m}});Q(\overrightarrow{X_{k}})\right) $ as%
\begin{eqnarray*}
	F_{n}(x):=\mathbb{S}_{n}\left( -x,-1;0,1\right) =\mathbb{Y}_{n-1}\left( -x,-1\right)
	=\sum_{n_{2}=0}^{\left[ \frac{n-1}{2}\right] }\binom{n-n_{2}-1}{n_{2}}%
	x^{n-2n_{2}-1}.
\end{eqnarray*}%
Let $m\in \mathbb{N}$. We define the following certain series of reciprocals
of the Fibonacci numbers:%
\begin{equation}
	\mathcal{R}_{m}(w)=\sum\limits_{j=1}^{\infty }\frac{w^{j}}{F_{mj}}.
	\label{ac5}
\end{equation}%
Let $\lambda +\mu =1$, $\lambda -\mu =\sqrt{5}$\ and $\lambda \mu =-1$.
Using (\ref{ac5}), we have%
\begin{eqnarray}
	\mathcal{R}_{m}(\lambda ^{2}w)-\mathcal{R}_{m}(\mu ^{2}w)
	=\sum\limits_{j=1}^{\infty }\frac{\lambda ^{2j}-\mu ^{2j}}{F_{mj}}w^{j} \label{YS2} 
	=\left( \lambda -\mu \right) \sum\limits_{j=1}^{\infty }\frac{F_{2j}}{%
		F_{mj}}w^{j}.
\end{eqnarray}%
Substituting $m=1$ into (\ref{ac5}), we get%
\begin{equation*}
	\mathcal{R}_{1}(w)=\sum\limits_{j=1}^{\infty }\frac{w^{j}}{F_{j}},
\end{equation*}%
(\textit{cf}. \cite{Gould}). 

Substituting $m=1$ into (\ref{YS2}), we also get 
\begin{equation*}
	\mathcal{R}_{1}(\lambda ^{2}w)-\mathcal{R}_{1}(\mu
	^{2}w)=\sum\limits_{j=1}^{\infty }\frac{F_{2j}}{F_{j}}w^{j}=G\left(
	w,-1,-1;1,2\right) .
\end{equation*}%
Combining the above equation with (\ref{ag-1}), we obtain%
\begin{equation*}
	\sum\limits_{j=1}^{\infty }\frac{F_{2j}}{F_{j}}w^{j}=\sum_{j=1}^{\infty }%
	\mathbb{S}_{j}\left( -1,-1;1,2\right) w^{j}.
\end{equation*}%
Comparing the coefficients of $w^{j}$ on both sides of the above equality,
we arrive at the following theorem:

\begin{theorem}
	Let $j\in \mathbb{N}$. Then we have%
	\begin{equation*}
		\mathbb{S}_{j}\left( -1,-1;1,2\right) =\frac{F_{2j}}{F_{j}}.
	\end{equation*}
\end{theorem}

Substituting $m=2$ and $w=1$ into (\ref{ac5}), the following certain series
of reciprocals of the Fibonacci numbers can be expressed in terms of the
Lambert series, which is defined by%
\begin{equation*}
	L(w)=\sum\limits_{j=1}^{\infty }\frac{w^{j}}{1-w^{j}}.
\end{equation*}%
Thus, we have
\begin{equation*}
	\sum\limits_{j=1}^{\infty }\frac{1}{F_{2j}}=L\left( \frac{3-\sqrt{5}}{2}%
	\right) -L\left( \frac{7-3\sqrt{5}}{2}\right) ,
\end{equation*}%
(\textit{cf}. \cite{Gould}).

\textbf{Lucas polynomials}: The Lucas polynomials can also be expressed in
terms of the polynomials $\mathbb{Y}_{n}\left( P(\overrightarrow{X_{m}}%
)\right) $ and $\mathbb{S}_{n}\left( P(\overrightarrow{X_{m}});Q(%
\overrightarrow{X_{k}})\right) $ as follows:
\begin{eqnarray*}
	L_{n}(x)=\mathbb{S}_{n}\left( -x,-1;2,-x\right) &=&2\mathbb{Y}_{n}\left( -x,-1\right)
	-x\mathbb{Y}_{n-1}\left( -x,-1\right) \\
	&=&2\sum_{n_{1}=0}^{\left[ \frac{n_{1}}{2}\right] }\binom{n_{1}-n_{2}}{n_{2}}%
	x^{n_{1}-2n_{2}}-\sum_{n_{1}=0}^{\left[ \frac{n_{1}-1}{2}\right] }\binom{%
		n_{1}-n_{2}-1}{n_{2}}x^{n_{1}-2n_{2}-2}. 
\end{eqnarray*}

Using (\ref{ac5}), we have%
\begin{eqnarray*}
	\mathcal{R}_{m}(\lambda ^{2}w)+\mathcal{R}_{m}(\mu ^{2}w)
	=\sum\limits_{j=1}^{\infty }\frac{\lambda ^{2j}+\mu ^{2j}}{F_{mj}}w^{j} 
	=\sum\limits_{j=1}^{\infty }\frac{L_{2j}}{F_{mj}}w^{j}.
\end{eqnarray*}%
Putting $m=4$ in the above equation after that combining the well-known
Binet's formulas for the Fibonacci numbers and the Lucas numbers, we get%
\begin{equation*}
	\mathcal{R}_{4}(\lambda ^{2}w)+\mathcal{R}_{4}(\mu ^{2}w)=\frac{1}{\lambda
		-\mu }\sum\limits_{j=1}^{\infty }\frac{1}{F_{2j}}w^{j}
\end{equation*}%
and
\begin{equation*}
	\sum\limits_{j=1}^{\infty }\frac{L_{2j}}{F_{4j}}w^{j}=\frac{1}{\lambda -\mu }%
	\sum\limits_{j=1}^{\infty }\frac{1}{F_{2j}}w^{j}.
\end{equation*}%
When $w=1$, the above equation reduces to the following result:

\begin{corollary}
	\begin{equation*}
		\sum\limits_{j=1}^{\infty }\frac{L_{2j}}{F_{4j}}=\frac{L\left( \frac{3-\sqrt{%
					5}}{2}\right) -L\left( \frac{7-3\sqrt{5}}{2}\right) }{\sqrt{5}}.
	\end{equation*}
\end{corollary}

\textbf{Fibonacci type polynomials and Lucas type polynomials:}

Fibonacci type polynomials and Lucas type polynomials can also be expressed
in terms of the polynomials $\mathbb{Y}_{n}\left( P(\overrightarrow{X_{m}}%
)\right) $ and $\mathbb{S}_{n}\left( P(\overrightarrow{X_{m}});Q(%
\overrightarrow{X_{k}})\right) $ as follows:

The Tribonacci-Lucas polynomials $t_{l,n}(x)$:
$$ t_{l,n}(x)=\mathbb{S}%
_{n}\left( -x^{2},-x,-1;3,-2x^{2},x\right) ,$$

the Tribonacci polynomials $T_{t,n}(x)$:
$$ T_{t,n}(x)=\mathbb{S}_{n}\left(
-x^{2},-x,-1;3,0,1\right)$$
(\textit{cf}. \cite{djordjevic,Rabolovic}). 

The Chebyshev polynomials of the first kind $T_{n}(x)$: $$T_{n}(x):=\mathbb{S}_{n}\left(
-2x,1;1,-x\right),$$

the Chebyshev polynomials of the second kind $U_{n}(x)$: $$U_{n}(x):=\mathbb{S}_{n}\left(
-2x,1;1\right) =\mathbb{Y}_{n}\left( -2x,1\right),$$

the Chebyshev polynomials of the third kind $T_{3,n}(x)$: $$T_{3,n}(x):=\mathbb{S}%
_{n}\left( -2x,1;1,-1\right),$$

the Chebyshev polynomials of the fourtht kind $T_{4,n}(x)$: $$T_{4,n}(x):=\mathbb{S}%
_{n}\left( -2x,1;1,1\right),$$

the monic $2$-orthogonal Chebyshev polynomial $\widehat{T}_{n}(x)$: $$\widehat{T}_{n}(x):=\mathbb{S}%
_{n}\left( x,\alpha ,\gamma \right),$$ where $\alpha $, $\gamma $ are
constants, for detail see \cite{alkanMTJp,Askey,Belbeshir,Bravo,AliB,Ceyvin,Charalambos,comtet,Deza,djordjevic,Dj-Mi2014,Forgacs0,Forgacs1,Forgacs2,gegenbauer,Goubi,GoubiMTJPAM,Gould,gould,Horadam,Hosonay,humbert,KahanBMS,KahanVit,KilarSimsekASCM,koshy1,KucukoAADm,KucukogluTJM,KucukogluAxio,Lia,lother,Munarini,Ortac,GulsahFilomat,GSG,Rabolovic,Pinch,Simsek,Simsekaadm,Simsek.mmas,SrivatavaChoi,srivastava,SrivatavaJNT,Suetin,Uygun}.

The Generalized Padovan sequences ($m$-Padovan numbers):%
\begin{equation*}
	\mathcal{P}_{n}^{(m)}:=\mathbb{S}_{n}\left( 0,\underset{m\text{ times }}{%
		\underbrace{-1,-1,\ldots ,-1}};0,1,1\right) 
\end{equation*}%
(\textit{cf}. \cite{Bravo}).

\textbf{Generating functions of the set }$F$\textbf{\ of words over }$%
\left\{ a,b\right\} $: The following generating functions of the set $%
\mathcal{F}$ of words over $\left\{ a,b\right\} $ without factor $w=a^{n}$
is given by%
\begin{eqnarray*}
	G\left( z,-2,\underset{m\text{ times}}{\underbrace{0,0,\ldots ,0}},1;x,%
	\underset{m-1\text{ times}}{\underbrace{0,0,\ldots ,0}},-1\right)
	&=&\sum\limits_{m=0}^{\infty }\mathbb{S}_{m}\left( -2,\underset{m\text{ times%
	}}{\underbrace{0,0,\ldots ,0}},1;x,\underset{m-1\text{ times}}{\underbrace{%
			0,0,\ldots ,0}},-1\right) z^{m} \\
	&=&\frac{1-z^{m}}{1-2z+z^{m+1}},
\end{eqnarray*}%
for detail, see \cite[p. 98]{lother}.

\subsection{Binet type formulas the polynomials $\mathbb{Y}_{n}\left( P(%
	\protect\overrightarrow{X_{m}})\right) $ and $\mathbb{S}_{n}\left( P(\protect%
	\overrightarrow{X_{m}});Q(\protect\overrightarrow{X_{k}})\right) $}

Here, we give Binet type formulas for the polynomials $\mathbb{Y}%
_{n}\left( P(\overrightarrow{X_{m}})\right) $ and $\mathbb{S}_{n}\left( P(%
\overrightarrow{X_{m}});Q(\overrightarrow{X_{k}})\right) $ with the aid of
generating functions.

Substituting $m=2$ into (\ref{Ag}), we now give Binet type formula for the
polynomials $\mathbb{Y}_{n}\left( P(\overrightarrow{X_{2}})\right) =\mathbb{Y%
}_{n}\left( P_{1}(x_{1}),P_{2}(x_{2})\right) $ as follows.

\begin{equation*}
	\frac{1}{1+P_{1}(x_{1})w+P_{2}(x_{2})w^{2}}=\sum_{n=0}^{\infty }\mathbb{Y}%
	_{n}\left( P_{1}(x_{1}),P_{2}(x_{2})\right) w^{n}.
\end{equation*}

The partial fraction decomposition of the left side of above equation is
given by%
\begin{eqnarray*}
	\sum_{n=0}^{\infty }\mathbb{Y}_{n}\left( P_{1}(x_{1}),P_{2}(x_{2})\right)
	w^{n} =\frac{1}{\left( 1-\frac{w}{a_{1}(x_{1},x_{2})}\right) \left( 1-%
		\frac{w}{a_{2}(x_{1},x_{2})}\right) },
\end{eqnarray*}%
where%
\begin{equation*}
	a_{1}(x_{1},x_{2})=\frac{-P_{1}(x_{1})+\sqrt{P_{1}^{2}(x_{1})-4P_{2}(x_{2})}%
	}{2P_{2}(x_{2})}
\end{equation*}%
and%
\begin{equation*}
	a_{2}(x_{1},x_{2})=\frac{-P_{1}(x_{1})-\sqrt{P_{1}^{2}(x_{1})-4P_{2}(x_{2})}%
	}{2P_{2}(x_{2})}.
\end{equation*}%
Therefore%
\begin{eqnarray*}
	\sum_{n=0}^{\infty }\mathbb{Y}_{n}\left( P_{1}(x_{1}),P_{2}(x_{2})\right)
	w^{n} &=&-\frac{a_{2}(x_{1},x_{2})}{a_{1}(x_{1},x_{2})-a_{2}(x_{1},x_{2})}%
	\sum_{n=0}^{\infty }\left( \frac{w}{a_{1}(x_{1},x_{2})}\right) ^{n} \\
	&&+\frac{a_{1}(x_{1},x_{2})}{a_{1}(x_{1},x_{2})-a_{2}(x_{1},x_{2})}%
	\sum_{n=0}^{\infty }\left( \frac{w}{a_{2}(x_{1},x_{2})}\right) ^{n}.
\end{eqnarray*}%
Comparing the coefficients of $w^{n}$ on both sides of the above equation,
we get%
\begin{eqnarray*}
	\mathbb{Y}_{n}\left( P_{1}(x_{1}),P_{2}(x_{2})\right) 
	=\frac{\left( a_{1}(x_{1},x_{2})\left( \frac{1}{a_{2}(x_{1},x_{2})}\right)
		^{n}-a_{2}(x_{1},x_{2})\left( \frac{1}{a_{1}(x_{1},x_{2})}\right)
		^{n}\right) }{a_{1}(x_{1},x_{2})-a_{2}(x_{1},x_{2})}.
\end{eqnarray*}%
After some elementary calculations in the above equation, the following
Binet type formula is obtained:

\begin{theorem}
	Let $n\in \mathbb{N}_{0}$. Then we have%
	\begin{eqnarray}
		&&\mathbb{Y}_{n}\left( P_{1}(x_{1}),P_{2}(x_{2})\right)   \label{BB1} \\
		&=&\frac{2^{n-1}\left( P_{2}(x_{2})\right) ^{n}}{\sqrt{%
				P_{1}^{2}(x_{1})-4P_{2}(x_{2})}}\left( \frac{P_{1}(x_{1})+\sqrt{%
				P_{1}^{2}(x_{1})-4P_{2}(x_{2})}}{\left( -P_{1}(x_{1})+\sqrt{%
				P_{1}^{2}(x_{1})-4P_{2}(x_{2})}\right) ^{n}}\right)   \notag \\
		&&-\frac{2^{n-1}\left( P_{2}(x_{2})\right) ^{n}}{\sqrt{%
				P_{1}^{2}(x_{1})-4P_{2}(x_{2})}}\left( \frac{P_{1}(x_{1})-\sqrt{%
				P_{1}^{2}(x_{1})-4P_{2}(x_{2})}}{\left( -P_{1}(x_{1})-\sqrt{%
				P_{1}^{2}(x_{1})-4P_{2}(x_{2})}\right) ^{n}}\right) .  \notag
	\end{eqnarray}
\end{theorem}

Substituting $m=2$ and $k=1$ into (\ref{ag-1}), we now give Binet type
formula for the polynomials $\mathbb{S}_{n}\left( P(\overrightarrow{X_{2}}%
);Q(\overrightarrow{X_{1}})\right) =\mathbb{S}_{n}\left(
P_{1}(x_{1}),P_{2}(x_{2});Q_{0}(x_{0}),Q_{1}(x_{1})\right) $ as follows:
\begin{equation*}
	\frac{Q_{0}(x_{0})+Q_{1}(x_{1})w}{1+P_{1}(x_{1})w+P_{2}(x_{2})w^{2}}%
	=\sum_{n=0}^{\infty }\mathbb{S}_{n}\left( P(\overrightarrow{X_{2}});Q(%
	\overrightarrow{X_{1}})\right) w^{n}.
\end{equation*}

The partial fraction decomposition of the left side of above equation is
given by%
\begin{eqnarray*}
	\sum_{n=0}^{\infty }\mathbb{S}_{n}\left( P(\overrightarrow{X_{2}});Q(%
	\overrightarrow{X_{1}})\right) w^{n} =\frac{Q_{0}(x_{0})+Q_{1}(x_{1})w}{%
		\left( 1-\frac{w}{a_{1}(x_{1},x_{2})}\right) \left( 1-\frac{w}{%
			a_{2}(x_{1},x_{2})}\right) },
\end{eqnarray*}%
where%
\begin{equation*}
	a_{1}(x_{1},x_{2})=\frac{-P_{1}(x_{1})+\sqrt{P_{1}^{2}(x_{1})-4P_{2}(x_{2})}%
	}{2P_{2}(x_{2})}
\end{equation*}%
and%
\begin{equation*}
	a_{2}(x_{1},x_{2})=\frac{-P_{1}(x_{1})-\sqrt{P_{1}^{2}(x_{1})-4P_{2}(x_{2})}%
	}{2P_{2}(x_{2})}.
\end{equation*}%
Therefore%
\begin{eqnarray*}
	&&\sum_{n=0}^{\infty }\mathbb{S}_{n}\left(
	P_{1}(x_{1}),P_{2}(x_{2});Q_{0}(x_{0}),Q_{1}(x_{1})\right) w^{n} \\
	&=&-\frac{a_{1}(x_{1},x_{2})a_{2}(x_{1},x_{2})\left( \frac{Q_{0}(x_{0})}{%
			a_{1}(x_{1},x_{2})}+Q_{1}(x_{1})\right) }{%
		a_{1}(x_{1},x_{2})-a_{2}(x_{1},x_{2})}\sum_{n=0}^{\infty }\left( \frac{1}{%
		a_{1}(x_{1},x_{2})}\right) ^{n}w^{n} \\
	&&+\frac{a_{2}(x_{1},x_{2})a_{1}(x_{1},x_{2})\left( \frac{Q_{0}(x_{0})}{%
			a_{2}(x_{1},x_{2})}+Q_{1}(x_{1})\right) }{\left(
		a_{1}(x_{1},x_{2})-a_{2}(x_{1},x_{2})\right) }\sum_{n=0}^{\infty }\left( 
	\frac{1}{a_{2}(x_{1},x_{2})}\right) ^{n}w^{n}.
\end{eqnarray*}%
Comparing the coefficients of $w^{n}$ on both sides of the above equation,
we get%
\begin{eqnarray*}
	&&\mathbb{S}_{n}\left(
	P_{1}(x_{1}),P_{2}(x_{2});Q_{0}(x_{0}),Q_{1}(x_{1})\right)  \\
	&=&-\frac{a_{1}(x_{1},x_{2})a_{2}(x_{1},x_{2})\left( \frac{Q_{0}(x_{0})}{%
			a_{1}(x_{1},x_{2})}+Q_{1}(x_{1})\right) }{%
		a_{1}(x_{1},x_{2})-a_{2}(x_{1},x_{2})}\left( \frac{1}{%
		a_{1}(x_{1},x_{2})}\right) ^{n} \\
	&&+\frac{a_{2}(x_{1},x_{2})a_{1}(x_{1},x_{2})\left( \frac{Q_{0}(x_{0})}{%
			a_{2}(x_{1},x_{2})}+Q_{1}(x_{1})\right) }{\left(
		a_{1}(x_{1},x_{2})-a_{2}(x_{1},x_{2})\right) }\left( 
	\frac{1}{a_{2}(x_{1},x_{2})}\right) ^{n}.
\end{eqnarray*}%
After some elementary calculations in the above equation, the following Binet type formula is obtained:

\begin{theorem}
	Let $n\in \mathbb{N}_{0}$. Then we have%
	\begin{eqnarray}
		&&\mathbb{S}_{n}\left(
		P_{1}(x_{1}),P_{2}(x_{2});Q_{0}(x_{0}),Q_{1}(x_{1})\right)   \label{BB2} \\
		&=&\frac{2^{n-1}\left( P_{2}(x_{2})\right) ^{n}}{\sqrt{%
				P_{1}^{2}(x_{1})-4P_{2}(x_{2})}}\left( \frac{\left( P_{1}(x_{1})+\sqrt{%
				P_{1}^{2}(x_{1})-4P_{2}(x_{2})}\right) Q_{0}(x_{0})-2Q_{1}(x_{1})}{\left(
			-P_{1}(x_{1})+\sqrt{P_{1}^{2}(x_{1})-4P_{2}(x_{2})}\right) ^{n}}\right)  
		\notag \\
		&&+\frac{2^{n-1}\left( P_{2}(x_{2})\right) ^{n}}{\sqrt{%
				P_{1}^{2}(x_{1})-4P_{2}(x_{2})}}\left( \frac{\left( -P_{1}(x_{1})+\sqrt{%
				P_{1}^{2}(x_{1})-4P_{2}(x_{2})}\right) Q_{0}(x_{0})+2Q_{1}(x_{1})}{\left(
			-P_{1}(x_{1})-\sqrt{P_{1}^{2}(x_{1})-4P_{2}(x_{2})}\right) ^{n}}\right) . 
		\notag
	\end{eqnarray}
\end{theorem}

Substituting $P_{1}(x_{1})=-1$, $P_{2}(x_{2})=-1$, $Q_{0}(x_{0})=0$ and $%
Q_{1}(x_{1})=1$ into (\ref{BB2}), we get the following explicit formula for the Fibonacci numbers:

\begin{corollary}
	Let $n\in \mathbb{N}_{0}$. Then we have%
	\begin{equation}
		\mathbb{S}_{n}\left( -1,-1;0,1\right) =\frac{(-2)^{n}}{\sqrt{5}}\left(\frac{1}{\left( 1-\sqrt{5}\right) ^{n}}- \frac{%
			1}{\left( 1+\sqrt{5}\right) ^{n}}%
		\right) =F_{n}.  \label{BF}
	\end{equation}
\end{corollary}

Substituting $P_{1}(x_{1})=-1$, $P_{2}(x_{2})=-1$, $Q_{0}(x_{0})=2$ and $%
Q_{1}(x_{1})=-1$ into (\ref{BB2}), we get the following explicit formula for the Lucas numbers:

\begin{corollary}
	Let $n\in \mathbb{N}_{0}$. Then we have%
	\begin{eqnarray}
		L_{n}=\mathbb{S}_{n}\left( -1,-1;2,-1\right) &=&2^n(-1)^{n}\left( 
		\frac{1}{\left( 1+\sqrt{5}\right) ^{n}}-\frac{1}{\left( 1-\sqrt{5}\right)
			^{n}}\right) .   \label{BL} 	 
	\end{eqnarray}
\end{corollary}
Substituting $P_{1}(x_{1})=-2$, $P_{2}(x_{2})=-1$, $Q_{0}(x_{0})=1$ and $%
Q_{1}(x_{1})=1$ into (\ref{BB2}),  we get the following corollary:

\begin{corollary}
	Let $n\in \mathbb{N}_{0}$. Then we have%
	\begin{eqnarray}
		\mathbb{S}_{n}\left( -2,-1;1,1\right) 
		=\frac{\left(-1\right)^n}{2\sqrt{2}} \left( \frac{-2+\sqrt{2}}{\left(1+\sqrt{2}\right)^{n}}+\frac{2+\sqrt{2}}{\left(1-\sqrt{2}\right)^{n}}\right).  \label{Bp}
	\end{eqnarray}
\end{corollary}

Observe that it is easy to see that the formula given in Eq. (\ref{Bp})
gives us the solution of Exercise 11 in Charalambides's book \cite[p. 267]%
{Charalambos}, given below:%
\begin{equation*}
	\mathbb{S}_{n}\left( -2,-1;1,1\right) =y_{n},
\end{equation*}%
where $y_{n}$ is the number of $n$-permutations of the set $\left\{
0,1,2\right\} $ with repetition and the restriction now two zeros and no two
ones are consecutive.

The Pell numbers%
\begin{equation*}
	\mathbb{S}_{n}\left( -2,-1;0,1\right) =P_{n}.
\end{equation*}

The Pell-Lucas numbers%
\begin{equation*}
	\mathbb{S}_{n}\left( -1,-1;2,-2\right) =Pl_{n},
\end{equation*}
(\textit{cf}. \cite{alkanMTJp}).

The Fibonacci polynomials%
\begin{equation*}
	\mathbb{S}_{n}\left( -x,-1;0,1\right) =F_{n}(x) .
\end{equation*}

The Lucas polynomials%
\begin{equation*}
	\mathbb{S}_{n}\left( -x,-1;2,-x\right) =L_{n}(x).
\end{equation*}%

\section{Higher order of the polynomials $\mathbb{Y}_{n}\left( P(\protect%
	\overrightarrow{X_{m}})\right) $ and $\mathbb{S}_{n}\left( P(\protect%
	\overrightarrow{X_{m}});Q(\protect\overrightarrow{X_{k}})\right) $} \label{section4}

In this section, we define higher order of the polynomials $\mathbb{Y}%
_{n}\left( P(\overrightarrow{X_{m}})\right) $ and $\mathbb{S}_{n}\left( P(%
\overrightarrow{X_{m}});Q(\overrightarrow{X_{k}})\right) $. We give some
properties of these polynomials with their special values.

The polynomials $\mathbb{Y}_{n} \left( P(%
\overrightarrow{X_{m}})\right) $ of order $\beta$, which is denoted by $\mathbb{Y}_{n}^{(\beta )}\left( P(%
\overrightarrow{X_{m}})\right) $, are defined by the following generating function
\begin{equation}
	F\left( w,P(\overrightarrow{X_{m}});\beta \right) =\frac{1}{\left(
		1+\sum\limits_{j=1}^{m}P_{j}(x_{j})w^{j}\right) ^{\beta }}%
	=\sum_{n=0}^{\infty }\mathbb{Y}_{n}^{(\beta )}\left( P(\overrightarrow{X_{m}}%
	)\right) w^{n}.  \label{Ag-h}
\end{equation}%
Here we note that%
\begin{equation*}
	\mathbb{Y}_{n}\left( P(\overrightarrow{X_{m}})\right) :=\mathbb{Y}%
	_{n}^{(1)}\left( P(\overrightarrow{X_{m}})\right) .
\end{equation*}

Some special values of the function $F\left( w,P(\overrightarrow{X_{m}});\beta
\right) $ are given as follows:%
\begin{equation*}
	F\left( w,-1;\beta \right) =\frac{1}{\left( 1-w\right) ^{\beta }}%
	=\sum_{n=0}^{\infty }\mathbb{Y}_{n}^{(\beta )}\left( -1\right) w^{n}.
\end{equation*}%
For $\left\vert w\right\vert <1$, using binomial theorem, we get%
\begin{equation*}
	\mathbb{Y}_{n}^{(\beta )}\left( -1\right) =(-1)^{n}\binom{-\beta }{n}=\binom{%
		\beta +n-1}{n}=\frac{\beta \left( \beta +1\right) \cdots \left( \beta
		+n-1\right) }{n!}=\frac{\left( \beta \right) _{n}}{n!},
\end{equation*}
where $\left( \beta \right) _{n}$ is denoted the rising factorial (the Pochhammer function). Few values of $\mathbb{Y}_{n}^{(\beta )}\left( -1\right) $ are given by%
\begin{eqnarray*}
	\mathbb{Y}_{n}^{(1)}\left( -1\right) &=&\frac{\left( 1\right) _{n}}{n!}=1, \\
	\mathbb{Y}_{n}^{(2)}\left( -1\right) &=&\frac{\left( 2\right) _{n}}{n!}=n+1,
	\\
	\mathbb{Y}_{n}^{(3)}\left( -1\right) &=&\frac{\left( 3\right) _{n}}{n!}=%
	\frac{\left( n+1\right) \left( n+2\right) }{2}, \\
	\mathbb{Y}_{n}^{(4)}\left( -1\right) &=&\frac{\left( 4\right) _{n}}{n!}=%
	\frac{\left( n+1\right) \left( n+2\right) (n+3)}{6} \\
	&&\vdots \\
	\mathbb{Y}_{n}^{(k)}\left( -1\right) &=&\frac{\left( k+1\right) _{n}}{n!}=%
	\binom{k+n}{n},
\end{eqnarray*}%
and so on.

Substituting $P_{1}(x)=-2x$, and $P_{2}(x)=1$, into (\ref{Ag-h}), we get the
most important family of orthogonal polynomials which are so-called the
ultraspherical polynomials or the Gegenbauer polynomials, for $0\leq
\left\vert x\right\vert <1$, $\beta >0$, 
\begin{equation*}
	C_{n}^{(\beta )}(x):=\mathbb{Y}_{n}^{(\beta )}\left( -2x,1\right) ,
\end{equation*}%
which are particular solutions of the Gegenbauer differential equation, and
these polynomials are also represented by the following Gaussian
hypergeometric series:%
\begin{equation*}
	\mathbb{Y}_{n}^{(\beta )}\left( -2x,1\right) =\frac{\left( 2\beta \right)
		_{n}}{n!}\text{ }_{2}F_{1}\left( -n,2\beta +n;\frac{1}{2}+\beta ;\frac{1-x}{2%
	}\right) ,
\end{equation*}%
and also special cases of the Jacobi polynomials are given by%
\begin{equation*}
	P_{n}^{\left( -\frac{1}{2}+\beta ,-\frac{1}{2}+\beta \right) }(x):=\frac{%
		\left( \frac{1}{2}+\beta \right) _{n}}{\left( 2\beta \right) _{n}}\mathbb{Y}%
	_{n}^{(\beta )}\left( -2x,1\right) .
\end{equation*}%
For detail about these
polynomials see also \cite{SrivatavaChoi}, and \cite{Suetin}.

With the aid of the Rodrigues formula for the Gegenbauer polynomials, we have%
\begin{eqnarray*}
	\mathbb{Y}_{n}^{(\beta )}\left( -2x,1\right) &=&\frac{(-1)^{n}\sqrt{\pi }%
		\Gamma (n+2\beta )(1-x^{2})^{-\beta +\frac{1}{2}}\binom{n-\frac{1}{2}}{n}}{%
		2^{n}\Gamma (2\beta )\Gamma \left( n+\beta +\frac{1}{2}\right) }\frac{d^{n}}{%
		dx^{n}}\left\{ (1-x^{2})^{n+\beta -\frac{1}{2}}\right\} \\
	&=&\frac{(-1)^{n}\sqrt{\pi }(n+1)C_{n}\Gamma (n+2\beta )(1-x^{2})^{-\beta +%
			\frac{1}{2}}}{2^{3n}\Gamma (2\beta )\Gamma \left( n+\beta +\frac{1}{2}%
		\right) }\frac{d^{n}}{dx^{n}}\left\{ (1-x^{2})^{n+\beta -\frac{1}{2}%
	}\right\} ,
\end{eqnarray*}%
where $C_{n}$ denotes the Catalan numbers, $\Gamma \left( \frac{1}{2}\right)
=\sqrt{\pi }$ and%
\begin{equation*}
	\Gamma \left( n+\frac{1}{2}\right) =\binom{n-\frac{1}{2}}{n}n!\sqrt{\pi }=%
	\frac{(2n)!}{4^{n}n!}\sqrt{\pi }.
\end{equation*}

When $\beta =\frac{1}{2}$ and $\beta =1$, the Gegenbauer polynomials reduce
to the Legendre polynomials and the Chebyshev polynomials of the second
kind, respectively.

The Rogers polynomials or the Rogers-Askey-Ismail polynomials, which are so
called continuous $q$-ultraspherical polynomials, are a family of orthogonal
polynomials introduced by Rogers (1892, 1893, 1894) in the course of his
work on the Rogers--Ramanujan identities, are defined in terms of the $q$%
-Pochhammer symbol and the basic hypergeometric series%
\begin{equation*}
	C_{n}^{(\beta \left\vert q\right. )}(x):=\mathbb{Y}_{n}^{(\beta \left\vert
		q\right. )}\left( -2x,1\right) =\frac{\left( \beta ;q\right) }{\left(
		q;q\right) }e^{in\theta }\text{ }_{2}\phi _{1}\left( q^{-n},\beta ;\beta
	^{-1}q^{1-n};q,q\beta ^{-1}e^{-2i\theta }\right) ,
\end{equation*}%
where $x=\cos \left( \theta \right) $,%
\begin{equation*}
	_{d}\phi _{c}\left( a_{1},\ldots ,a_{d};b_{1},\ldots ,b_{c};q,z\right) =\sum 
	\frac{(a_{1},\ldots ,a_{d};q)_{n}}{(b_{1},\ldots ,b_{c};q)_{n}}\left(
	(-1)^{n}q^{\binom{n}{2}}\right) ^{1+c-d}z^{n}
\end{equation*}%
and%
\begin{equation*}
	(a_{1},\ldots ,a_{d};q)_{n}=(a_{1};q)_{n}(a_{2};q)_{n}\cdots
	(a_{d};q)_{n},(a_{1};q)_{n}
\end{equation*}%
denotes the $q$-Pochhammer symbol ($q$-shifted factorial) or $q$-analog of
the Pochhammer symbol $(x)_{n}$, 
\begin{equation*}
	(a_{j};q)_{n}=\left( 1-a_{j}\right) \left( 1-a_{j}q\right) \cdots \left(
	1-a_{j}q^{n-1}\right)
\end{equation*}%
and $(a_{j};q)_{0}=1$ and also%
\begin{equation*}
	\lim_{q\rightarrow 1}\frac{(q^{x};q)_{n}}{(1-q)^{n}}=(x)_{n}
\end{equation*}%
(\textit{cf}. \cite{Askey}).

Substituting $\beta =\frac{1}{2}$, $P_{1}(x)=-2x$ and $P_{2}(x)=1$ into (\ref%
{Ag-h}), we get other most important family of orthogonal polynomials which
are so-called the Legendre polynomials%
\begin{equation*}
	P_{n}(x):=\mathbb{Y}_{n}^{(\frac{1}{2})}\left( -2x,1\right) .
\end{equation*}

Substituting $\beta =-\frac{1}{2}$, $P_{1}(x)=-3x,P_{2}(x)=0$, and $%
P_{3}(1)=1$ into (\ref{Ag-h}), we get the Pincherle polynomials introduced
by Humbert in 1921,%
\begin{equation*}
	P_{n}(x):=\mathbb{Y}_{n}^{(-\frac{1}{2})}\left( -3x,0,1\right)
\end{equation*}%
and%
\begin{equation*}
	\mathbb{Y}_{2n}^{(-\frac{1}{2})}\left( -3x,0,1\right) =\frac{d^{n}}{dx^{n}}%
	\left\{ x^{n}(x^{2}-1)^{n}\right\}
\end{equation*}%
(\textit{cf}. \cite{humbert}).

Substituting $P_{1}(x)=-mx,P_{2}(x)=P_{m}(x)=0,\ldots ,P_{m}(1)=1$ into (\ref%
{Ag-h}), we get the Humbert polynomials, which are generalization of the
Pincherle polynomials introduced by Humbert in 1921,%
\begin{equation*}
	\Pi _{n,m}^{\beta }(x):=\mathbb{Y}_{n}^{(\beta )}\left( -mx,\underset{m-1%
		\text{ times}}{\underbrace{0,0,\ldots ,0}},1\right) .
\end{equation*}%
For other important properties and applications of these polynomials, the
references given here, among other references, may be reviewed (\textit{cf}. \cite{Askey,Charalambos,comtet,Dj-Mi2014,SrivatavaChoi,srivastava,Suetin}).

The Fibonacci type polynomials of higher order in two variables $\mathcal{G}%
_{v}^{(h)}(x,y,k,m,n)$:%
\begin{equation*}
	\mathcal{G}_{v}^{(h)}(x_{1},x_{m+n},k,m,n):=\mathbb{Y}_{v}^{(h)}\left(
	-x_{1}^{k},\underset{m+n-1\text{ times}}{\underbrace{0,0,\ldots ,0}}%
	,-x_{m+n}^{m}\right)
\end{equation*}%
(\textit{cf}. \cite{GulsahFilomat}).

By using (\ref{Ag-h}), we get%
\begin{equation*}
	\mathbb{Y}_{n}^{(\beta +\gamma )}\left( P(\overrightarrow{X_{m}})\right)
	=\sum\limits_{j=0}^{n}\mathbb{Y}_{n}^{(\beta )}\left( P(\overrightarrow{X_{m}%
	})\right) \mathbb{Y}_{n-j}^{(\gamma )}\left( P(\overrightarrow{X_{m}}%
	)\right) .
\end{equation*}

The polynomials $\mathbb{S}%
_{n}\left( P(\overrightarrow{X_{m}});Q(\overrightarrow{%
	X_{k}})\right) $  of order $\alpha,\beta$, which is denoted by $\mathbb{S}%
_{n}^{(\alpha ,\beta )}\left( P(\overrightarrow{X_{m}});Q(\overrightarrow{%
	X_{k}})\right) $ , are defined by the following generating function
\begin{eqnarray}
	H\left( w,P(\overrightarrow{X_{m}});Q(\overrightarrow{X_{k}});\alpha
	,\beta \right)  \label{ag-1y} &=&\frac{\left( \sum\limits_{j=0}^{k}Q_{j}(x_{j})w^{j}\right) ^{\alpha }}{%
		\left( 1+\sum\limits_{j=1}^{m}P_{j}(x_{j})w^{j}\right) ^{\beta }}\\
	&=&\sum_{n=0}^{\infty }\mathbb{S}_{n}^{(\alpha ,\beta )}\left( P(%
	\overrightarrow{X_{m}});Q(\overrightarrow{X_{k}})\right) w^{n}.  \notag
\end{eqnarray}%
It is clear that%
\begin{equation*}
	\mathbb{S}_{n}\left( P(\overrightarrow{X_{m}});Q(\overrightarrow{X_{k}}%
	)\right) :=\mathbb{S}_{n}^{(1,1)}\left( P(\overrightarrow{X_{m}});Q(%
	\overrightarrow{X_{k}})\right) .
\end{equation*}

Combining (\ref{ag-1y}) with (\ref{Ag-h}), we get the following functional
equation:%
\begin{equation*}
	H\left( w,P(\overrightarrow{X_{m}});Q(\overrightarrow{X_{k}});1,\beta
	\right) =\left( \sum\limits_{j=0}^{k}Q_{j}(x_{j})w^{j}\right) F\left( w,P(%
	\overrightarrow{X_{m}});\beta \right) .
\end{equation*}

Combining (\ref{ag-1y}) with (\ref{ag-1}), we get%
\begin{equation*}
	\sum_{n=0}^{\infty }\mathbb{S}_{n}^{(1,\beta )}\left( P(\overrightarrow{X_{m}%
	});Q(\overrightarrow{X_{k}})\right) w^{n}=\left(
	\sum\limits_{j=0}^{k}Q_{j}(x_{j})w^{j}\right) \sum_{n=0}^{\infty }\mathbb{Y}%
	_{n}^{(\beta )}\left( P(\overrightarrow{X_{m}})\right) w^{n}.
\end{equation*}%
Therefore%
\begin{equation*}
	\sum_{n=0}^{\infty }\mathbb{S}_{n}^{(1,\beta )}\left( P(\overrightarrow{X_{m}%
	});Q(\overrightarrow{X_{k}})\right) w^{n}=\sum\limits_{j=0}^{k}\sum_{n=j}^{\infty
	}Q_{j}(x_{j})\mathbb{Y}_{n}^{(\beta )}\left( P(\overrightarrow{X_{m}}%
	)\right) w^{n}.
\end{equation*}%
Comparing the coefficients of $w^{n}$ on both sides of the above
equality, we arrive at the following theorem:

\begin{theorem}
	Let $n,k\in \mathbb{N}$ with $k\leq n$. Then we have%
	\begin{equation}
		\mathbb{S}_{n}^{(1,\beta )}\left( P(\overrightarrow{X_{m}});Q(%
		\overrightarrow{X_{k}})\right) =\sum\limits_{j=0}^{k}Q_{j}(x_{j})\mathbb{Y}%
		_{n-j}^{(\beta )}\left( P(\overrightarrow{X_{m}})\right) .  \label{Ag2}
	\end{equation}%
\end{theorem}
The generalized Catalan polynomials $\mathcal{P}_{v,m}^{h,Q_{1}}(x_{1})$:%
\begin{equation*}
	\mathcal{P}_{v,m}^{h,Q_{1}}(x_{1}):=\mathbb{S}_{v}^{(1,h)}\left( -m,\underset%
	{m-1\text{ times}}{\underbrace{0,0,\ldots ,0}},-x_{1};1,Q_{1}(x_{1})\right)
\end{equation*}%
(\textit{cf}. \cite{Goubi}).

\textbf{Generating functions for two }$2$\textbf{-variable Simsek polynomials:%
} Replacing $\left( \sum\limits_{j=0}^{k}Q_{j}(x_{j})w^{j}\right) ^{\alpha }$
and $\left( 1+\sum\limits_{j=1}^{m}P_{j}(x_{j})w^{j}\right) ^{\beta }$ by $%
\left( 1+\lambda w\right) ^{\alpha _{1}}\left( 1+\delta w^{2}\right)
^{\alpha _{2}}$ and $\lambda -1+\lambda ^{2}w$ in the equation (\ref{ag-1y}%
), repectively, Khan et al. (\cite{KahanBMS}, \cite{KahanVit}) defined the
following generating function:%
\begin{equation*}
	\frac{\left( 1+\lambda w\right) ^{\alpha _{1}}\left( 1+\delta w^{2}\right)
		^{\alpha _{2}}}{\lambda -1+\lambda ^{2}w}=\sum_{n=0}^{\infty }\mathbb{S}%
	_{n}^{(\alpha _{1}+\alpha _{2},1)}\left( \lambda -1,\lambda ^{2};1,1,\delta
	\right) w^{n},
\end{equation*}%
where%
\begin{equation*}
	Y_{n}(\alpha _{1},\alpha _{2};\lambda ,\delta ):=\mathbb{S}_{n}^{(\alpha
		_{1}+\alpha _{2},1)}\left( \lambda -1,\lambda ^{2};1,1,\delta \right)
\end{equation*}%
denotes the two $2$-variable Simsek polynomials and $Y_{n}(0,0;\lambda
,\delta )$ denotes the Simsek numbers. $Y_{n}(\alpha ;\lambda
):=Y_{n}(\alpha _{1},0;\lambda ,\delta )$ known as the Simsek polynomials (%
\textit{cf}. for detail, see \cite{KucukoAADm,KucukogluTJM,KucukogluAxio,Simsek,Simsekaadm,SrivatavaJNT}).

\section{Recurrence relation of the polynomials $\mathbb{Y}	_{n}\left( P(\protect\overrightarrow{X_{m}})\right) $ and $\mathbb{S}_{n}\left( P(\protect\overrightarrow{X_{m}});Q(\protect\overrightarrow{X_{k}})\right) $} \label{section5}

In this section, using partial derivative equations of the generating functions for the polynomials  $\mathbb{Y}%
_{n}\left( P(\overrightarrow{X_{m}})\right) $ and $\mathbb{S}_{n}\left( P(%
\overrightarrow{X_{m}});Q(\overrightarrow{X_{k}})\right)$, we give recurrence relations of these polynomials.

Some differential equations of the functions $F\left( w,P(\overrightarrow{%
	X_{m}})\right) $ and  $G\left( w,P(\overrightarrow{X_{m}});Q(\overrightarrow{%
	X_{k}})\right) $ are given as follows:%
\begin{equation}
	\frac{\partial }{\partial w}\left\{ F\left( w,P(\overrightarrow{X_{m}}%
	)\right) \right\} =-\sum\limits_{j=1}^{m}jP_{j}(x_{j})w^{j-1}F\left( w,P(%
	\overrightarrow{X_{m}});2\right) ,  \label{ad-1}
\end{equation}%
\begin{eqnarray}
	\frac{\partial }{\partial w}\left\{ G\left( w,P(\overrightarrow{X_{m}});Q(%
	\overrightarrow{X_{k}})\right) \right\}
	&=&\sum\limits_{l=1}^{k}lQ_{l}(x_{l})w^{l-1}F\left( w,P(\overrightarrow{X_{m}%
	})\right)  \label{ad0} \\
	&&-\sum\limits_{l=0}^{k}Q_{l}(x_{l})w^{l}\frac{\partial }{\partial w}\left\{
	F\left( w,P(\overrightarrow{X_{m}})\right) \right\} ,  \notag
\end{eqnarray}%
and%
\begin{eqnarray}
	&&F\left( w,P(\overrightarrow{X_{m}});2\right) \frac{\partial }{\partial w}%
	\left\{ G\left( w,P(\overrightarrow{X_{m}});Q(\overrightarrow{X_{k}})\right)
	\right\}  \label{ad2} \\
	&=&\sum\limits_{l=1}^{k}lQ_{l}(x_{l})w^{l-1}+\sum\limits_{l=1}^{k}\sum%
	\limits_{j=1}^{m}lQ_{l}(x_{l})P_{j}(x_{j})w^{j+l-1}-\sum\limits_{l=0}^{k}%
	\sum\limits_{j=1}^{m}jQ_{l}(x_{l})P_{j}(x_{j})w^{l+j-1}.  \notag
\end{eqnarray}

The above differential equation can also be given by%
\begin{eqnarray}
	\frac{\partial }{\partial w}\left\{ G\left( w,P(\overrightarrow{X_{m}});Q(%
	\overrightarrow{X_{k}})\right) \right\}
	&=&\sum\limits_{l=1}^{k}lQ_{l}(x_{l})w^{l-1}F\left( w,P(\overrightarrow{X_{m}%
	})\right) -  \label{ad3} \\
	&&-\sum\limits_{l=0}^{k}\sum%
	\limits_{j=1}^{m}jQ_{l}(x_{l})P_{j}(x_{j})w^{l+j-1}F\left( w,P(%
	\overrightarrow{X_{m}});2\right) .  \notag
\end{eqnarray}%
or%
\begin{eqnarray}
	\frac{\partial }{\partial w}\left\{ G\left( w,P(\overrightarrow{X_{m}});Q(%
	\overrightarrow{X_{k}})\right) \right\}
	&=&\sum\limits_{l=1}^{k}lQ_{l}(x_{l})w^{l-1}F\left( w,P(\overrightarrow{X_{m}%
	})\right)  \label{ad4} \\
	&&-\sum\limits_{j=1}^{m}jP_{j}(x_{j})w^{j-1}G\left( w,P(\overrightarrow{X_{m}%
	});Q(\overrightarrow{X_{k}})\right) F\left( w,P(\overrightarrow{X_{m}}%
	)\right) .  \notag
\end{eqnarray}

Combining (\ref{ad-1}) with (\ref{Ag}) and (\ref{Ag-h}), we get%
\begin{equation*}
	\sum_{n=1}^{\infty }n\mathbb{Y}_{n}\left( P(\overrightarrow{X_{m}})\right)
	w^{n-1}=-\sum\limits_{j=1}^{m}jP_{j}(x_{j})\sum_{n=0}^{\infty }\mathbb{Y}%
	_{n}^{(2)}\left( P(\overrightarrow{X_{m}})\right) w^{n+j-1}.
\end{equation*}%
Therefore%
\begin{equation*}
	\sum_{n=0}^{\infty }(n+1)\mathbb{Y}_{n+1}\left( P(\overrightarrow{X_{m}}%
	)\right) w^{n}=-\sum\limits_{j=1}^{m}jP_{j}(x_{j})\sum_{n=j-1}^{\infty }%
	\mathbb{Y}_{n-j+1}^{(2)}\left( P(\overrightarrow{X_{m}})\right) w^{n}.
\end{equation*}%
After the necessary operations in the previous equation, the $w^{n}$
coefficients on both sides of this equation are compared and the following
result is easily obtained:

\begin{theorem}
	Let $n\in \mathbb{N}_{0}$. Then we have%
	\begin{equation*}
		\mathbb{Y}_{n+1}\left( P(\overrightarrow{X_{m}})\right) =-\frac{1}{n+1}%
		\sum\limits_{j=1}^{m}jP_{j}(x_{j})\mathbb{Y}_{n-j+1}^{(2)}\left( P(%
		\overrightarrow{X_{m}})\right) .
	\end{equation*}
\end{theorem}

Combining (\ref{ad0})with (\ref{Ag}) and (\ref{ag-1}), we get%
\begin{eqnarray*}
	\sum_{n=1}^{\infty }n\mathbb{S}_{n}\left( P(\overrightarrow{X_{m}});Q(%
	\overrightarrow{X_{k}})\right) w^{n-1}
	&=&\sum\limits_{l=1}^{k}lQ_{l}(x_{l})\sum_{n=0}^{\infty }\mathbb{Y}%
	_{n}\left( P(\overrightarrow{X_{m}})\right) w^{n+l-1}\\
	&&-\sum\limits_{l=0}^{k}Q_{l}(x_{l})\sum_{n=1}^{\infty }n\mathbb{Y}%
	_{n}\left( P(\overrightarrow{X_{m}})\right) w^{n+l-1}.
\end{eqnarray*}%
After the necessary operations in the previous equation, the $w^{n}$
coefficients on both sides of this equation are compared and the following
result is easily obtained:

\begin{theorem}
	Let $n\in \mathbb{N}$. Then we have%
	\begin{eqnarray*}
		\mathbb{S}_{n+1}\left( P(\overrightarrow{X_{m}});Q(\overrightarrow{X_{k}}%
		)\right) &=&\frac{1}{n+1}\sum\limits_{l=1}^{k}lQ_{l}(x_{l})\mathbb{Y}%
		_{n+1-l}\left( P(\overrightarrow{X_{m}})\right) \\
		&&-\sum\limits_{l=0}^{k}\left( \frac{n+1-l}{n+1}\right) Q_{l}(x_{l})\mathbb{Y%
		}_{n-l+1}\left( P(\overrightarrow{X_{m}})\right) .
	\end{eqnarray*}
\end{theorem}

Note that if the above operations are applied to the derivative equations
from Eq. (\ref{ad2}) to Eq. (\ref{ad4}), other different
recurrence relations are obtained. We omit these relations and solution
extraction here.



\end{document}